\documentclass[11pt]{article}
\usepackage[greek,english]{babel}
\usepackage[iso-8859-7]{inputenc}
\usepackage{amssymb}
\usepackage{amsmath}
\usepackage{amsthm}
\usepackage{latexsym}
\usepackage{amsfonts}
\usepackage{graphicx}
\usepackage{graphics}

\theoremstyle{plain}
\theoremstyle{definition}

\newcommand{\bbb}[1]{\mbox{\boldmath$#1$}}

\newcommand{\OO} {{\varOmega}}
\newcommand{\oo} {{\omega}}

\newcommand{\bi} {{\beta}}

\newcommand{\Ga} {{\varGamma}}

\newcommand{\ld} {{\ldots}}
\newcommand{\sm} {{\smallsetminus}}
\newcommand{\thi} {{\theta}}

\newcommand{\de} {{\delta}}

\newcommand{\Si} {{\varSigma}}
\newcommand{\la} {{\lambda}}
\newcommand{\el} {{\ell}}

\newcommand{\vPi} {{\varPi}}

\newcommand{\f} {{\varphi}}
\newcommand{\mi} {{\mu}}

\newcommand{\dis}{\displaystyle}
\newcommand{\ssum}{\sum\limits}
\newcommand{\ct}{{\cal{T}}}
\newcommand{\cu}{{\cal{U}}}
\newcommand{\cp}{{\cal{P}}}

\newcommand{\ch}{{\cal{H}}}
\newcommand{\cm}{{\cal{M}}}

\newcommand{\ra}{{\rightarrow}}
\newcommand{\oD}{{\overline{D}}}
\newcommand{\oA}{{\overline{A}}}
\newcommand{\oB}{{\overline{B}}}

\newcommand{\fa}{{\forall}}

\def\1{\it1\hspace*{-0.150cm}{\footnotesize{I}}}

\def\R{{\mathbb{R}}}
\def\C{{\mathbb{C}}}

\def\N{{\mathbb{N}}}

\begin{document}
\title{\bf Common hypercyclic vectors and\\ universal functions}
\author{\bf G. Costakis - N. Tsirivas}\footnotetext{\hspace{-0.5cm}}
\footnotetext{{The research project is implemented within the framework of the Action ``Supporting Postdoctoral Researchers'' of the Operational Program ``Educational and Lifelong Learning'' (Action's Beneficiary: General Secretariat for Research and Technology), and is co-financed by the European Social Fund (ESF) and the Greek State.}}
\date{}
\maketitle
\noindent
{\bf Abstract:} Let $X,Y$ be two separable Banach or Fr\'{e}ch\`{e}t spaces, and $T_n:X\ra Y$ be a sequence from linear and continuous operators. We say that the sequence $(T_n)$, $n=1,2,\ld$ is universal, if there exists some vector $v\in X$ such that the sequence $T_n(v)$, $n=1,2,\ld$ is dense in $Y$. If $X=Y$ we say that the sequence $(T_n)$ is hypercyclic.

More generally we consider an uncountable subset $A$ from complex numbers and for every fixed $a\in A$ we consider a sequence $(T_{a,n})$, $n=1,2,\ld$, from linear and continuous operators, $T_{a,n}:X\ra Y$.

The problem of common universal vectors is whether the uncountable family of sequences of operators $(T_{a,n})$, $n=1,2,\ld$ for $a\in A$ share a common universal vector.

We examine, in this work, some specific cases of this problem for translation, differential and backward shift operators. We study also some approximating problems about universal Taylor series.
\vspace*{0.2cm} 
\vspace*{0.2cm} 
\section{Introduction}\label{sec1}
\noindent

The basic aim of the present paper is to develop briefly the results of our research project with the number: PE1 4126, K.A. 3486, that is a research project of the action ``Supporting Postdoctoral researchers'' of the General Cecretariat for Research and Technology and is cofinanced by the European Social Fund (ESF) and the Greek State. This paper is a formal obligation of our postdoctoral research project.

As for the results of our project we have to remark the following:

We have made progress in all the components in which we are obligated by the regulations, to do.

More specifically: 

In the first component we give a positive answer in an open question of G. Costakis from 2007, and we strictly extent the results of G. Costakis from his paper \cite{6}, in our paper \cite{21}.

In the components 2 and 3 the results are optimal with the sense that we have found a characterization of the sequences with which we have a positive result for a non degenerate subinterval of the set of complex numbers $\C$. It is quite surprizing that the condition is the same for the two problems and when the condition holds, then automatically we have a positive result for a subset of full 2-dimensional Lebesque measure. 

In components 4 and 5 we deal with some approximating problems that come from universal Taylor series. The most important paper and simultaneously the starting point for the area of universal Taylor series is the paper of Vassili Nestoridis \cite{19}. Another important paper in the area is \cite{18} a wealthy in results paper.

In our component 4 we deal with an approximating problem in the complex plane concerning weighted universal Taylor series see \cite{25}. 

In an other supplementary paper \cite{26} for this problem with give the best possible, in a specific sense, lower bound of the well known in the literature asymptotic convergence factor. This number is appeared, with a natural way, in Krylov's subspace method for solving large-sparse systems of linear equations, that is 1 one from the top 10 numerical methods with many applications in natural sciences as mechanics, oceanografy, meteorology and so on. The bibliography and the references for this subject is vast. We only refer a paper \cite{9} that shows the above iteraction.

We are able to give two different proofs for this lower bound. One using potential theory \cite{26} and one other using the theory of universal Taylor series \cite{25}. So, we are happy that we have found applications of universal Taylor series by giving such an inequality.

{\bf In our component 5 we prove an approximating result about polynomials. More specifically: Let $\bbb{D}$ be the open unit disc. Let $(\la_n)_{n\in\N}$ be a strictly increasing subsequence from natural numbers. We fix a strictly increasing subsequence from natural numbers, $\bbb{(m_n)_{n\in\N}}$. Then there exists a sequence of polynomials $\bbb{p_n(z)}$ of the form $\bbb{p_n(z)=z^{m_n}q_n(z)}$, $\bbb{n=1,2,\ld}$,
where $\bbb{q_n(z)}$, $\bbb{n\in\N}$ is a sequence from complex polynomials,
with $\bbb{deg(q_n(z))=\la_n}$ for $\bbb{n\in\N}$ such that the sequence $\bbb{(p_n)}$, $\bbb{n\in\N}$ is dense in $\bbb{A(K)}$ for every compact set $\bbb{K\subset D^c}$, with connected complement if and only if $\bbb{\underset{n\ra+\infty}{\lim\sup}\Big(\dfrac{\la_n}{m_n}\Big)=+\infty}$.}

Using this approximating theorem we prove the existence of doubly universal Taylor series \cite{8}.

As for the methods of this project. We have used knowledge from some areas as Functional Analysis, Complex analysis, measure theory, multivariable Real analysis, topology, but above all potential theory. We have used potential theory for the three out of the seven papers. Potential theory is a powerful tool for solving quite a lot problems from universal Taylor series.

Significant work in this area has made by Stephen Gardiner and other researchers see \cite{10}, \cite{11}, \cite{12}, \cite{13} and \cite{14}.

As a conclusion: The most of our results are the best possible in a sense.

As an overall we have written 7 papers in our project. One of them has already published and the other six have already published in math arXiv. Three of them \cite{21}, \cite{22} and \cite{24} have already taken positive reports for publication.

We begin with the respective definitions and terminology.

Let $X$ be some non-void set and $T:X\ra X$ be a map.

If we take some $y\in X$ we define the orbit of $y$ under $T$ be the following sequence from elements in $X$,
\[
T(y),T(T(y)),\ld\underbrace{T(T(\ld(}_{n-\text{times}}T(y))\ld),\ld\;.
\]

Typically we define the ``iterates'' of map $T$ as follows:
\[
T^1:=T \ \ \text{and} \ \ T^{n+1}:=T\circ T^n \ \ \text{for} \ \ n=1,2,\ld,
\]
where with $T\circ T^n$ we mean the usual composition of maps, $T$ and $T^n$.

The orbit of $y$ under the sequence of maps $(T^n)$, $n=1,2,\ld$, is the set
\[
A_{T,y}:=\{\oo\in X\mid\exists\; n\in\N:\;\oo=T^n(y)\}.
\]
Of course $A_{T,y}\subset X$ for every $y\in X$. More generally we define the orbits as follows:

Let $X,Y$ be some non-void sets and $T_i:X\ra Y$, for $i\in I$ be a family of maps where $I$ is a set of indices, finite, denumerable or uncountable.

The orbit $O(y,(T_i)_{i\in I})$ of some $y\in X$, under the family of maps $(T_i)_{i\in I}$ is defined to be the set
\[
O(y,(T_i)_{i\in I}):=\{\oo\in Y\mid\exists\;i\in I:\;\oo=T_i(y)\}.
\]
Now, we suppose that the set $Y$ has a topology $\ct_Y$.

We say that the family of maps $(T_i)_{i\in I}$ is universal if there exists some element $y\in X$ such that:
\[
\overline{O(y,(T_i)_{i\in I})}=Y
\]
where the closure is taken with respect to the topology $\ct_Y$. Then the element $y$ is called universal under the family $(T_i)_{i\in I}$. For the above definitions see \cite{3}, \cite{15} and \cite{16}.

We denote
\[
\cu((T_i)_{i\in I})=\{y\in X\mid\;\overline{O(y,(T_i)_{i\in I})}=Y\},
\]
the set of universal elements of $X$ under the action of the family $(T_i)_{i\in I}$.
We can express the definition of universality as an equation also as follows. We fix the family $(T_i)_{i\in I}$, $T_i:X\ra Y$, $i\in I$.

We consider the map $F:X\ra\cp(Y)$ where with $\cp(Y)$ we denote the powerset of $Y$ with the formula:
\[
F(y):=\overline{O(y,(T_i)_{i\in I})}, \ \ \text{for} \ \ y\in X.
\]
The element $y$ is universal under the family $(T_i)_{i\in I}$ whether $y$ is a solution of the equation:
\[
F(y)=Y.  \eqno{\mbox{$(\ast)$}}
\]
By the above we consider the set
\[
L((T_i)_{i\in I}):=\{y\in X\mid\;F(y)=Y\}
\]
of solutions of the equation $(\ast)$, that is equal with the set of universal elements of $X$ under $(T_i)_{i\in I}$, that is $\cu((T_i)_{i\in I})=L((T_i)_{i\in I})$. The first example of universality seems to go back to M. Fekete in 1914 (quoted in \cite{20}) who discovered the existence of a universal Taylor series $\ssum^{+\infty}_{n=1}a_nt^n$: for any continuous function $g$ on $[-1,1]$ with $g(0)=0$, there exists an increasing sequence of integers $(m_k)$ such that $\ssum^{m_k}_{n=1}a_nt^n\ra g(t)$ uniformly as $k\ra+\infty$. Here $X=\C^\N$, $Y$ is the space of all continuous functions on $[-1,1]$ vanishing at 0, and
\[
T_i((a_n))=\sum^i_{n=1}a_nt^n, \ \ i=1,2,\ld\;.
\]
Since then, universal families have been exhibited in a huge number of situations see \cite{15}.

The most interesting examples of universality are concrete, as usually in mathematics. Below we denote $\C$, the set of complex numbers, $\N=\{1,2,\ld\}$ the set of natural numbers, $D:=\{z\in\C\mid|z|<1\}$ the open unit disc,\\
$\oD:=\{z\in\C\mid|z|\le1\}$ the closed unit disc, $A^c:=\C\sm A$ the complement of a subset $A$ of $\C$, $\overset{\circ}{A}$ the interior of a subset $A$ of $\C$, $\oA$ the closure of a subset $A$ in some topological space, $A(K):=\{f:K\ra\C\mid f$ is continuous and holomorphic in $\overset{\circ}{K}\}$ for some compact subset $K$ of $\C$. We consider the set $A(K)$ endowed with the supremum norm with which is become a Banach algebra.

Usually the sets $X,Y$ are some well known Fr\'{e}ch\`{e}t or Banach spaces with obvious topologies and the family $(T_i)_{i\in I}$ is a sequence from some well known linear and continuous operators. In the case where $X=Y$ the universal family $(T_i)_{i\in I}$ is called hypercyclic and every universal element is called hypercyclic, also.

In this case we write $\ch C((T_n)_{n\in\N})$ instead of $\cu((T_i)_{i\in I})$.

In every one from the paragraphs of this paper we describe the result of our project. At the end we give only the necessary references for this exposition. Further references there are in the references, of the papers \cite{21}-\cite{26}.
\section{Common hypercyclic functions for translation operators}\label{sec2}
\noindent

In this paragraph we describe the results of Problem 1 - Component 1 of our project. The first result for hypercyclicity of translation operators goes back to a result of G. Birkhoff \cite{5}.

We denote $\ch(\C)$ the set of entire functions.

We consider the set $\ch(\C)$ endowed with the topology of local uniform convergence $\ct_u$. The space $(\ch(\C),\ct_u)$ is a Fr\'{e}ch\`{e}t space so Baire's Category Theorem holds in this space. We remind that a set in a topological space is $G_\de$ if it can be written as a denumerable intersection of open sets, and a set is dense if its closure is whole the space.

For every fixed complex number $a$ we consider the translation entire function\\
$t_a:\C\ra\C$ with the formula $t_a(z)=z+a$ for every $z\in\C$. We fix some $a\in\C$.

The translation operator $T_a:\ch(\C)\ra\ch(\C)$ with the formula:
\[
T_a(f):=f\circ t_a \ \ \text{for every} \ \ f\in\ch(\C),
\]
is a well defined linear and continuous operator, where with $f\circ t_a$ we denote the usual composition of the functions $f$ and $t_a$. We symbolize $T^1_a:=T_a$ and $T^{n+1}_a:=T_a\circ T^n_a$ for $n\in\N$, where with $T_a\circ T^n_a$, $n\in\N$ we denote the usual composition of the operators $T_a$ and $T^n_a$. It holds $T^n_a=T_{na}$ for every $a\in\C\sm\{0\}$ and $n\in\N$. It is well known that the sequence $(T^n_a)$, $n\in\N$ is hypercyclic, when $a\neq0$.

More generally we consider a sequence from complex numbers $(a_n)$, such that $a_n\ra\infty$ as $n\ra\infty$. Then it is well known that the sequence $(T_{a_n})$, $n\in\N$ is hypercyclic \cite{17}. It is also known that the set $\ch C((T_{a_n})_{n\in\N}))$ is $G_\de$ and dense. Let $(b_m)$, $m=1,2$, be a sequence from non-zero complex numbers. Then from this result the set $\ch C((T_{b_ma_n})_{n\in\N})$ is $G_\de$ and dense for every $m\in\N$. So, by Baire's Category Theorem the set $\bigcap\limits^\infty_{m=1}\ch C((T_{b_ma_n})_{n\in\N})$ is $G_\de$ and dense. However, if we have an uncountable subset $I$ of $\C\sm\{0\}$ we do not know whether the intersection $\bigcap\limits_{b\in I}\ch C((T_{ba_n})_{n\in\N})$ is non-void because Baire's Category Theorem holds only for a countable intersection of open dense sets.

G. Costakis proved in 2007 the following result: see \cite{6}.

We consider a sequence $(\la_n)_{n\in\N}$ of complex numbers that satisfies the following condition $(\Si)$. For every $M>0$ there exists a subsequence $(\mi_n)_{n\in\N}$ of $(\la_n)_{n\in\N}$ such that
\begin{enumerate}
\item[(i)] $|\mi_{n+1}|-|\mi_n|>M$ for every $n\in\N$ and
\item[(ii)] $\ssum^{+\infty}_{k=1}\dfrac{1}{|\mi_k|}=+\infty$.
\end{enumerate}
We denote $C(0,1):=\{z\in\C\mid|z|=1\}$. Then the intersection
\[
\bigcap_{b\in C(0,1)}\ch C((T_{b\la_n})_{n\in\N})
\]
is a residual set, that is, it contains a $G_\de$, dense set.

{\bf Costakis asked in the same paper \cite{6} whether
\[
\bbb{\dis\bigcap_{b\in\C\sm\{0\}}\ch C((T_{b\la_n})_{n\in\N})\neq\emptyset.}
\]
We proved in \cite{21} that this question has positive reply.}

F. Bayart examined similar results in $\R^n$, see \cite{2}.

{\bf Further we consider a sequence $\bbb{(\la_n)_{n\in\N}}$ from complex numbers  that satisfies the following condition $\bbb{(C')}$. For every positive number $\bbb{M>0}$, there exists a subsequence $\bbb{(\mi_n)_{n\in\N}}$ of $\bbb{(\la_n)_{n\in\N}}$ such that:
\begin{enumerate}
\item[(i)] $\mi_1\neq0$.
\item[(ii)] $\bbb{|\mi_{n+1}|-|\mi_n|>M}$ for every $\bbb{n\in\N}$ and
\item[(iii)] $\bbb{\underset{n\ra+\infty}{\lim\sup}|\mi_n|\bigg(\ssum^{+\infty}_{k=n}\dfrac{1}
    {|\mi_k|}\bigg)=+\infty}$.
\end{enumerate}
In our paper \cite{22} we proved that
\[
\bbb{\dis\bigcap_{b\in C(0,1)}\ch C((T_{b\la_n})_{n\in\N})}
\]
is residual, hence non-empty.}

This result is of course a strict extention of the result in \cite{6} because it holds for sequences $(\la_n)$ with polynomial growth, for example the sequence $\la_n=n^2$ $n\in\N$.
\section{Common hypercyclic vectors for certain families of differential operators}\label{sec3}
\noindent

In this paragraph we describe the results of Problem 2 - Component 2 of our project.

Now, we consider the Fr\'{e}ch\`{e}t space $\ch(\C)$ of entire functions endowed with the topology of local uniform convergence.

Let $\la\in\C$ fixed. We consider the dilation function $\f_\la:\C\ra\C$ with the formula:
\[
\f_\la(z):=\la z \ \ \text{for every} \ \ z\in\C.
\]
Of course $\f_\la$ is an entire function. We consider the $n$-th derivative operator inductively with the formula.

$D^1:\ch(\C)\ra\ch(\C)$, where $D^1(f)=f'$ where $f\in\ch(\C)$ and $f'$ is the usual derivative of $f$
\[
D^{n+1}:=D^1\circ D^n, \ \ D^n:\ch(\C)\ra\ch(\C) \ \ \text{for every} \ \ n\in\N \ \ \text{where} \ \ D^1\circ D^n
\]
be the usual composition of the operators $D^1$ and $D^n$ for $n\in\N$.

That is $D^n(f)=f^{(n)}$ where $f^{(n)}$ is the usual $n$-th derivative of $f\in\ch(\C)$.

Now, for every fixed $\la\in\C$ we consider the derivative operators $(T_{\la,n})$, $n\in\N$ $T_{\la,n}:\ch(\C)\ra\ch(\C)$, with the formula $T_{\la,n}(f)=D^n(f
\circ\f_\la)$ for every $n\in\N$, $f\in\ch(\C)$, where with $f\circ\f_\la$ we mean the composition of entire functions $f$ and $\f_\la$.

It is well known that for every fixed $\la\in\C\sm\{0\}$, the sequence $(T_{\la,n})_{n\in\N}$ is hypercyclic and the set $\ch C((T_{\la,n})_{n\in\N})$ is $G_\de$, dense.\\
Costakis-Sambarino proved \cite{7} that the set $\bigcap\limits_{\la\in\C\sm\{0\}}\ch C((T_{\la,n})_{n\in\N})$ is residual.

It is also known that for every fixed $\la\in\C\sm\{0\}$ and $(m_n)_{n\in\N}$ a strictly increasing subsequence of natural numbers the sequence $(T_{\la,m_n})_{n\in\N}$ is hypercyclic and the set $\ch C((T_{\la,m_n})_{n\in\N})$ is $G_\de$, dense. Let $I\subset\C\sm\{0\}$ and $I$ be uncountable. Is it true that
\[
\bigcap_{\la\in I}\ch C((T_{\la,m_n})_{n\in\N})\neq\emptyset\,?
\]
This question concerns Problem 2 - Component 2 of our project.

As for this problem we have proved in \cite{23} the following result: Let $(m_n)_{n\in\N}$ be a strictly increasing subsequence from natural numbers.

{\bf 1) When $\bbb{(m_n)_{n\in\N}}$ is a sequence such that the series $\bbb{\ssum^{+\infty}_{k=1}\dfrac{1}{m_n}}$ converges to a positive number, then we have that
\[
\bbb{\dis\bigcap_{\la\in[\thi_1,\thi_2]}\ch C((T_{\la,m_n})_{n\in\N})=\emptyset,}
\]
for every interval $\bbb{[\thi_1,\thi_2]}$ where $\bbb{0<\thi_1<\thi_2<+\infty}$.

2) When $\bbb{(m_n)}$ is a sequence such that $\bbb{\ssum^{+\infty}_{n=1}\dfrac{1}{m_n}=+\infty}$ then we have\\ $\bbb{\bigcap\limits_{\la\in(0,+\infty)}\ch C((T_{\la,m_n})_{n\in\N})}$ is a residual set.

3) When $\bbb{(m_n)_{n\in\N}}$ is a sequence such that $\bbb{\ssum^{+\infty}_{n=1}\dfrac{1}{m_n}=+\infty}$ then there exists a full measure subset $\bbb{I\subset\C\sm\{0\}}$ such that $\bbb{\la_2(\C\sm I)=0}$, and
\[
\bbb{\dis\bigcap_{\la\in I}\ch C((T_{\la,m_n})_{n\in\N})\neq\emptyset,}
\]
where with $\bbb{\la_2}$ we mean the 2-dimensional Lebesque measure on the complex plane.}

Our method here is based on the specific features of our problem. Another possible method that we have not tried is to apply a general common hypercyclicity criterion. See \cite{4} for some general powerful criteria, alternatively.
\section{Common hypercyclic vectors for families of Backward shift operators}\label{sec4}
\noindent

In this paragraph we describe the results of Problem 3 - Component 3 of our project.

Now, we consider the Banach space $\el^2$ of square summable sequences over the field of complex numbers $\C$, endowed with the topology that is induced by the $\el^2$ norm $\|\cdot\|_2:\el^2\ra\R^+$, where
\[
\|x\|_2:=\bigg(\sum^{+\infty}_{j=1}|x_j|^2\bigg)^{1/2} \ \ \text{for every} \ \ x=(x_1,x_2,\ld)\in\el^2.
\]
We write $\|\cdot\|:=\|\cdot\|_2$ for simplicity. We remark that here we have a Banach space whereas in the two previous Problems 1 and 2 we have a Fr\'{e}ch\`{e}t space.

Let $B$ be the unweighted backward shift operator on $\el^2$, that is
\[
B((x_1,x_2,x_3,\ld))=(x_2,x_3,\ld), \ \ \text{for every} \ \ (x_1,x_2,\ld)\in\el^2.
\]
For every fixed $\la\in\C$, $|\la|>1$ we consider the sequence of operators $(\la B)^n:\el^2\ra\el^2$, $n\in\N$, that is defined inductively as follows:
\begin{align*}
&(\la B)^1:=\la B \ \ \text{and} \\
&(\la B)^{n+1}:=(\la B)\circ(\la B)^n \ \ \text{for every} \ \ n\in\N, \ \ \text{where with} \ \ (\la B)\circ(\la B)^n
\end{align*}
we mean the usual composition of the operators $(\la B)$ and $(\la B)^n$ for every $n\in\N$.

It is a well known result that the sequence $(\la B)^n$ $n\in\N$ is hypercyclic for every $\la\in\C$, $|\la|>1$.

Abakumov and Gordon \cite{1} proved the best result about common hypercyclicity for these operators, that is they proved
\[
\bigcap_{|\la|>1}\ch C((\la B)^n\;n\in\N)\neq\emptyset.
\]
Later on, Costakis and Sambarino \cite{7} gave a different proof of this result, which, roughly speaking, is based on the so called common hypercyclicity criterion. In this criterion, Baire's category theorem appears. Actually, Costakis and Sambarino showed that $\bigcap\limits_{|\la|>1}\ch C((\la B)^n\; n\in\N)$ is a $G_\de$ and dense subset of $(\el^2,\|\cdot\|)$; hence non-empty. What is interesting here is the uncountable range of $\la$'s, which makes things harder if one wishes to apply Baire's theorem.

Our Problem 3 is a refinement of the previous result.

Let $(k_n)_{n\in\N}$ be a fixed strictly increasing subsequence of natural numbers. It is known, and very easy to prove, that the sequence $(\la B)^{k_n}$ $n\in\N$ is also hypercyclic, that is there exists $x\in\el^2$ such that the set $((\la B)^{k_n}(x)$ $n\in\N)$ is dense in $\el^2$ for every fixed $\la$, $|\la|>1$. It holds that the set $\ch C((\la B)^{k_n}\;n\in\N)$ is a $G_\de$ and dense subset of $(\el^2,\|\cdot\|)$ for every fixed $\la$, $|\la|>1$.

In our Problem 3 we have posed the following question:\\
Fix a strictly increasing subsequence $(k_n)_{n\in\N}$ of natural numbers. For which uncountable sets $J\subset\{\la\in\C\mid|\la|>1\}$,
\[
\bigcap_{\la\in J}\ch C((\la B)^{k_n}\;n\in\N)\neq\emptyset\,?
\]
Our main result in \cite{24} is the following:\vspace*{0.2cm} \\
\noindent
{\bf Theorem.}
{\bf Let $\bbb{(k_n)_{n\in\N}}$ be a strictly increasing subsequence of natural numbers.
\begin{enumerate}
\item[(i)] If $\bbb{\ssum^{+\infty}_{n=1}\dfrac{1}{k_n}<+\infty}$ then $\bbb{\bigcap\limits_{\la\in I}\ch C((\la B)^{k_n}\;n\in\N)=\emptyset}$ for every non-degenerate interval $\bbb{I\subset\C\sm\oD}$.
\item[(ii)] If $\bbb{\ssum^{+\infty}_{n=1}\dfrac{1}{k_n}=+\infty}$ then the set $\bbb{\bigcap\limits_{\la\in(1,+\infty)}\ch C((\la B)^{k_n}\;n\in\N)}$ is residual in $\bbb{\el^2}$; hence
 \[
\bbb{\dis\bigcap_{\la\in(1,+\infty)}\ch C((\la B)^{k_n}\;n\in\N)\neq\emptyset.}
\]
\item[(iii)] If $\bbb{\ssum^{+\infty}_{n=1}\dfrac{1}{k_n}=+\infty}$, there exists a $\bbb{G_\de}$ and dense subset $\bbb{\cp}$ in $\bbb{\{\la\in\C\mid|\la|>1\}}$ with full 2-dimensional Lebesque measure in $\bbb{\{\la\in\C\mid|\la|>1\}}$ such that
\[
\bbb{\dis\bigcap_{\la\in\cp}\ch C((\la B)^{k_n}\;n\in\N)}
\]
is residual in $\bbb{\el^2}$. In particular, $\bbb{\bigcap\limits_{\la\in\cp}\ch C((\la B)^{k_n}\;n\in\N)\neq\emptyset}$, where
\[
\bbb{\la_2((\C\sm\oD)\sm\cp)=0},
\]
\end{enumerate}
}
\section{Universal Taylor serious on specific compact sets and a general lower bound for the asymptotic convergence factor}\label{sec5}
\noindent

In this paragraph we desicribe the results of Problem 5 - Component 4 of our project.

In the first three problems of our project we have the problem to find common hypercyclic vectors for an uncountable set of sequences of operators for which we know, in advance, that every one from the sequences of operators is hypercyclic.

In this problem we differentiate from the previous three problems. We consider the more general sense of universal families of operators and not hypercyclic operators. We consider only one sequence of linear and continuous operators between two specific spaces for which we want to prove that this is universal.

More specifically:

We consider the set $\ch(D)$ of holomorphic functions on $D$, that is,\\
$\ch(D):=\{f:D\ra\C\mid f$ is holomorphic$\}$.

Let some $f\in\ch(D)$. We denote $S_n(f)$, $n\in\N$, the $n$-th partial sum of the Taylor's development of $f$ about 0, that is
\[
S_n(f)(z):=\sum^n_{k=0}\frac{f^{(k)}(0)}{k!}z^k,  \ \ \text{for every} \ \ z\in\C, \ \ n\in\N.
\]
Let $\cp(\C)$ be the powerset of $\C$ and
\[
\cm_{D^c}:=\{K\in\cp(\C)\mid K \ \ \text{is compact}, \ \ K^c \ \ \text{is connected}, \ \ K\cap D=\emptyset\}.
\]

Let some $K\in\cm_{D^c}$. We consider the set $A(K)$ endowed with the supremum norm.

Without loss of generality we write below $S_n(f)=S_n(f)\upharpoonright K$ for some $f\in\ch(D)$ for every $n\in\N$, without distinction for simplicity, where $S_n(f)\upharpoonright K$ is the restriction of $S_n(f)$ on $K$.

In his very influential paper \cite{19}, Vassili Nestoridis proved that there exists some $f\in\ch(D)$, such that for every $K\in\cm_{D^c}$ we have that:
\[
\overline{\{S_n(f),\,n\in\N\}}=A(K).
\]
We can write this result in the language of universality as follows:

The set $\cm_{D^c}$ has the cardinality of continuum as the closed interval $[0,1]$.

Let $j:[0,1]\ra\cm_{D^c}$ be a map $1-1$ and onto.

So we have $\cm_{D^c}:=\{j(\thi)\mid\thi\in[0,1]\}$ or we can write simply that $\cm_{D^c}$ is the family $j(\thi)$, $\thi\in[0,1]$.

For every fixed $\thi\in[0,1]$, we consider the sequence of operators
\[
T_{\thi,n}:\ch(D)\ra A(j(\thi)), \ \ n\in\N
\]
that is defined by the formula:
\[
T_{\thi,n}(f)(z):=S_n(f)\upharpoonright_{j(\thi)}(z)=S_n(f)(z) \ \ \text{for every} \ \ n\in\N, \ \ f\in\ch(D), \ \ z\in j(\thi).
\]

Let $\cu((T_{\thi,n})\;n\in\N)$ be the set of universal vectors of the sequence of operators $T_{\thi,n}$, $n\in\N$. That is,
\[
\cu((T_{\thi,n})\;n\in\N):=\{f\in\ch(D)\mid\overline{\{T_{\thi,n}(f),n\in\N\}}=A(j(\thi))\}.
\]
We denote
\[
\cu(D):=\bigcap_{\thi\in[0,1]}\cu((T_{\thi,n})\;n\in\N).
\]
It is proved in \cite{19} that $\cu(D)\neq\emptyset$. We consider the space $\ch(D)$ endowed with the topology of local uniform convergence. Even if this result seems to be a problem of common universality as our results in Problems 1, 2 and 3, this is not true because as it is proved in \cite{19} every set $\cu((T_{\thi,n})\;n\in\N)$ is $G_\de$, dense in $\ch(D)$ for every $\thi\in[0,1]$ and also there exists some sequence of numbers $(\thi_n)_{n\in\N}$ in $[0,1]$ such that:
\[
\bigcap_{\thi\in[0,1]}\cu((T_{\thi,n})\;n\in\N)=\bigcap^{+\infty}_{n=1}\cu((T_{\thi_n,m})\;m\in\N),
\]
so this differentiate this problem from the problems of common universality, because in order to have that $\cu(D)$ is non-empty it suffices only to have that every one from the sets $\cu((T_{\thi,n})\;n\in\N)$ is non void for every $\thi\in[0,1]$ whereas something similar does not hold in the first three problems as we have shown, that implies that the first three problems are strictly common hypercyclicity problems.

Let some sequence $(\bi_n)_{n\in\N}$ from complex numbers. Let some $K\in\cm_{D^c}$, $K=j(\thi)$, $\thi\in[0,1]$. We consider the sequence of operators $T_{(\bi\nu)},\thi,n:\ch(D)\ra A(K)$, for $n\in\N$, that is defined with the formula:
\[
T_{(\bi\nu)},\thi,n(f)(z):=\bi_nS_n(f)(z), \ \ n\in\N, \ \ f\in\ch(D), \ \ z\in j(\thi).
\]
Let $\cu((T_{(\bi\nu)},\thi,n)\;n\in\N)$ be the set of universal vectors for the sequence\\ $(T_{(\bi\nu)}\;\thi,n)$ $n\in\N$. We set
\[
\cu(\bi\nu):=\bigcap_{\thi\in[0,1]}\cu((T_{(\bi\nu)},\thi,n)_{n\in\N}).
\]
There exists a sequence $(\thi_n)$ of numbers in $[0,1]$ such that:
\[
\bigcap_{\thi\in[0,1]}\cu((T_{(\bi\nu)},\thi,n)_{n\in\N})=\bigcap^{+\infty}_{m=1}
\cu((T_{(\bi\nu)},\thi_m,n)_{n\in\N}).
\]
In \cite{27} it is proved that $\cu(\bi\nu)$ is non-empty if and only if the sequence $(|\bi n|^{1/n})_{n\in\N}$ has 1 as a limit point, and this holds if and only if every one from the sets $\cu((T_{(\bi\nu)},\thi,n)_{n\in\N})$ is non-empty.

So, if we take a sequence $(\bi_n)$ such that the sequence $(|\bi_n|^{1/n})_{n\in\N}$ does not have 1 as a limit point then $\cu(\bi_n)=\emptyset$. By the above equivalent conditions we get that there exists some $\thi_0\in[0,1]$ (at least):
\[
\cu((T_{(\bi_\nu)},\thi_0,n)_{n\in\N})=\emptyset.
\]
So the following question is arising naturally.

Let some sequence $(\bi_n)_{n\in\N}$ such that the sequence $(|\bi_n|^{1/n})_{n\in\N}$ does not have 1 as a limit point.

Does there exists some positive number $\thi_1\in[0,1]$ such that
\[
\cu((T_{(\bi_n)},\thi_1,n)_{n\in\N})\neq\emptyset\,?
\]
In our paper \cite{25} we examine the above problem and we formulate now the results of this paper.

Let $K\in\cm_{D^c}$ such that $K=\bigcup\limits^{n_0}_{i=0}K_i$, $n_0\in\N$, and $K_i$ be the connected components of $K$, $K_0=\oD$, and every one from the components contains more than one point. Let some subset $A\subseteq\C$ and $f:A\ra\C$ be a complex function. We denote
\[
\|f\|_A:=\sup\{x\in\R\mid\exists\,a\in A:x=|f(a)|\}\in[0,+\infty].
\]
We denote $V_n$, the set of complex polynomials of degree at most $n$ for $n\in\N$.
We take some complex polynomials $p_j$, $j=0,1,\ld,n_0$ different each other and we consider the function
\[
F:K\ra\C \ \ \text{where} \ \ F(z)=p_j(z) \ \ \text{if} \ \ z\in K_j \ \ \text{for every} \ \ j=0,1,\ld,n.
\]
We denote:
\[
d(V_n,F):=\min\{\|F-p\|_K,p\in V_n\} \ \ \text{for} \ \ n\in\N.
\]
It is well known that the number $\underset{n\ra+\infty}{\lim\sup}d(V_n,F)^{1/n}\in(0,1)$ and it is a number independent from the function $F$ and it is depending only on the compact set $K$.

We symbolize $\rho_K:=\underset{n\ra+\infty}{\lim\sup}d(V_n,F)^{1/n}$.

We fix a sequence $(\bi_n)_{n\in\N}$ from complex numbers. We set $\vPi:=K\sm K_0=\bigcup\limits^{n_0}_{i=1}K_i$. Let $\thi_0\in[0,1]$ be the unique number such that: $j(\thi_0)=\vPi$. We consider the space $A(\oD)$ endowed with the supremum norm, that is a Banach algebra. We consider the sequence of linear and continuous operators:
\[
T_{(\bi_\nu),\thi_0,n}:A(K_0)\ra A(\vPi) \ \ \text{for} \ \ n\in\N
\]
that is defined by the formula:
\[
T_{(\bi_\nu),\thi_0,n}(f):\vPi\ra\C
\]
for every $f\in A(K_0)$, $n\in\N$, where
\[
T_{(\bi_\nu),\thi_0,n}(f)(z):=\bi_nS_n(f)(z)
\]
for every $n\in\N$, $f\in A(K_0)$, $z\in\vPi=j(\thi_0)$.

Let $\cu((T_{(\bi_\nu),\thi_0,n})_{n\in\N})$ be the set of universal vectors for the sequence of operators $(T_{(\bi_\nu),\thi_0,n})_{n\in\N}$.
\vspace*{0.2cm} \\
\noindent
{\bf Theorem 1)} {\bf If the sequence $\bbb{(|\bi_n|^{1/n})_{n\in\N}}$ has a limit point in the interval $\bbb{\Big(\rho_K,\dfrac{1}{\rho_K}\Big)}$, then the set $\bbb{\cu((T_{(\bi_n),\thi_0,n})_{n\in\N})}$ is $\bbb{G_\de}$, dense.

Let $\bbb{h:=dist(\{0\},\vPi)}$}.\vspace*{0.2cm} \\
\noindent
{\bf Theorem 2)} {\bf If $\bbb{\underset{n\ra+\infty}{\lim\sup}|\bi_n|^{1/n}<\dfrac{1}{h}}$ then
\[
\bbb{\cu((T_{(\bi_n),\thi_0,n})_{n\in\N})=\emptyset.}
\]
Let $\bbb{\OO:=(\C\sm\vPi)\cup\{\infty\}}$. We consider the set $\bbb{\OO}$ endowed with the one point compactification of the usual topology of $\bbb{\C}$ and the unique Green's function of $\bbb{\OO}$ with pole at infinity $\bbb{g_\OO}$.

We set
\[
\bbb{M_K:=e^{\dis\max_{z\in\oD}g_\OO(z,\infty)}.}
\]
}
\noindent
{\bf Theorem 3)} {\bf If $\bbb{\underset{n\ra+\infty}{\lim\inf}|\bi_n|^{1/n}>M_K}$, then
\[
\bbb{\cu((T_{(\bi_n),\thi_0,n})_{n\in\N})=\emptyset.}
\]
In addition to the above results in \cite{26} we have found a significant inequality that is
\[
\bbb{\rho_L\ge\max\limits_{0\le i\le n_0}\sup\limits_{z\in\overset{\circ}{K}_i}\dfrac{dist(z,K^c_i)}{dist(z,K\sm K_i)}} \eqno{\mbox{$\bbb{(\ast)}$}}
\]
that gives us a lower bound for the very well known in the literature number $\rho_L$, that is optimal in a sense.}
\section{Doubly universal Taylor series}\label{sec6}
\noindent

In this paragraph we describe the results of Problem 6 - Component 5 of our project.

This work was born by our effort to solve the following problem.

Let $K$ be some compact subset of $\C$ with connected complement.

By the well known Mergelyan's Theorem we have that the set of all complex polynomials is dense in $A(K)$ when this is endowed with the supremum norm $\|\cdot\|_\infty$. The set $V:=\{p\mid p$ is a complex polynomial$\}$ of all complex polynomials is a very big subset of $A(K)$. The problem is to find strictly smaller than $V$ subsets of $V$ that remain dense in the space $(A(K),\|\cdot\|_\infty)$ also.

Trying to solve this problem we proved in \cite{8} the following Theorem. We fix two strictly increasing subsequences $(m_n)_{n\in\N}$ and $(\la_n)_{n\in\N}$ from natural numbers such that $\underset{n\ra+\infty}{\lim\sup}\dfrac{\la_n}{m_n}=+\infty$.

Then there exists a sequence of polynomials $(p_n)$, $n\in\N$ such that %
\[
p_n(z)=z^{m_n}q_n(z)\ \ \text{for every}\ \ n\in\N, \ \ z\in\C,
\]
where $q_n(z)$ is a sequence of polynomials such that $deg(q_n(z))=\la_n$ for every $n\in\N$ such that the sequence $(p_n)$ is dense in $(A(K),\|\cdot\|_\infty)$, for every $K\in\cm_{D^c}$.

Using the previous approximating result we solved a result of universality as follows.

Below we describe our result using the terminology that we have developed in the previous Section \ref{sec5}.

Let $(X,\ct_X)$, $(Y,\ct_Y)$ be two topological spaces. We consider the set $X\times Y$ endowed with the product topology $\ct_{X\times Y}$, that is inherited by the topologies $\ct_X$ and $\ct_Y$. Let $A\subseteq X$ and $B\subset Y$. Then it is well known that
\begin{eqnarray}
\oA\times\oB=\overline{A\times B}  \label{eq1}
\end{eqnarray}
where the closure in the second member of this equality is taken with respect to the product topology. Such an equality holds for finite, denumerable or non-denumerable number of topological spaces, as it is well known.

Using the previous equality we can take a result of universality on the product space $A(K_1)\times A(K_2)$, where $K_1,K_2\in\cm_{D^c}$.

More specifically:

Let some $f\in\cu(D)$ where the set $\cu(D)$ is the set of universal Taylor series on the unit disc in Nestoridis paper \cite{19} as we have defined it in page 11 of our present paper.

Let some $K_1,K_2\in\cm_{D^c}$.

We set
\[
B(K_1,f):=\{S_n(f):K_1\ra\C,\;n\in\N\}
\]
\[
B(K_2,f):=\{S_m(f):K_2\ra\C,\;m\in\N\}.
\]
We have:
\[
B(K_1,f)\times B(K_2,f):=\{(S_n(f),S_m(f))\mid
\]
\[
S_n(f):K_1\ra\C,\;S_m(f):K_2\ra\C,\;n,\;m\in\N\}.
\]
Because of $f\in\cu(D)$, we have $\overline{B(K_1,f)}=A(K_1)$, $\overline{B(K_2,f)}=A(K_2)$ and by the previous equality (\ref{eq1}) we get
$\overline{B(K_1,f)\times B(K_2,f)}=A(K_1)\times A(K_2)$,
that means that from the universality on the spaces $A(K_1)$ and $A(K_2)$ we take immediately a universality result on the space\\ $A(K_1)\times A(K_2)$. A question that is arised naturally is the following:

Can we have some strictly subset $\Ga\subseteq B(K_1,f)\times B(K_2,f)$ such that\\ $\overline{\Ga}=A(K_1)\times A(K_2)$.

Or in other words, can we have universality with a strict smaller set than $B(K_1,f)\times B(K_2,f)$\,?

It is obvious that the candidate set $\Ga$ have to be infinite.

We solved this problem for $\Ga(f,(\la_n)_{n\in\N}):=\{(S_n(f),S_{\la_n}(f)\mid n\in\N\}$ for some specific $f\in\cu(D)$, where $(\la_n)_{n\in\N}$ is a strictly increasing subsequence of natural numbers.

More specifically we proved in our published paper \cite{8} the following result:

We fix a strictly increasing subsequence $(\la_n)_{n\in\N}$ of natural numbers.

Let
\[
\cu(D,(\la_n)_{n\in\N}):=\{f\in\ch(D)\mid\;\fa\,K\in\cm_{D^c}:\overline{\Ga(f,(\la_n)_{n\in\N})}
=A(K)^2\}.
\]
Then we have \cite{8}.\vspace*{0.2cm}\\
\noindent
{\bf Theorem.}  {\bf The set $\bbb{\cu(D,(\la_n)_{n\in\N})}$ is non-empty if and only if $\bbb{\underset{n\ra+\infty}{\lim\sup}(\la_n/n)=+\infty}$. If $\bbb{\cu(D,(\la_n)_{n\in\N})\neq\emptyset}$ then this set is $\bbb{G_\de}$ and dense in $\bbb{(\ch(D),\ct_u)}$, where with $\bbb{\ct_u}$ we denote the topology of local uniform convergence on the space $\bbb{\ch(D)}$}.

We note that Vagia Vlachou solved \cite{28} the above problem in the general case for arbitrary finite number of spaces.

\vspace*{1cm}
\noindent
{\sc{University of Crete,\\
Department of Mathematics and Applied Mathematics,\\
GR-700 13 Heraklion,\\
Crete,\\ Greece}}\bigskip\\
e-mail addresses: \\
\noindent
costakis@uoc.gr, \;\;tsirivas@uoc.gr

\end{document}